\documentclass{amsart}
\usepackage{amssymb}
\usepackage{graphicx}
\usepackage{amscd}
\usepackage{amsmath}
\usepackage{amsfonts}

\newtheorem{theorem}{Theorem}
\theoremstyle{plain}

\newtheorem{corollary}{Corollary}

\newtheorem{definition}{Definition}

\newtheorem{lemma}{Lemma}

\newtheorem{proposition}{Proposition}
\newtheorem{remark}{Remark}

\numberwithin{equation}{section}

\input{tcilatex}

\begin{document}
\author{William Chin and Leonid Krop}
\address{DePaul University\\
Chicago, Illinois 60614\\
USA}
\email{{\small wchin@condor.depaul.edu , lkrop@condor.depaul.edu }}
\date{This work was supported in part by a grant from the University
Research Council of DePaul University.}
\date{}
\date{This work was supported in part by a grant from the University
Research Council of DePaul University.}
\date{}
\date{This work was supported in part by a grant from the University
Research Council of DePaul University.}
\date{}
\date{This work was supported in part by a grant from the University
Research Council of DePaul University.}
\date{}
\title[Spectra of quantized hyperalgebras]{Spectra of quantized hyperalgebras%
}
\date{This work was supported in part by a grant from the University
Research Council of DePaul University.}
\date{}
\date{This work was supported in part by a grant from the University
Research Council of DePaul University.}
\date{}
\date{This work was supported in part by a grant from the University
Research Council of DePaul University.}
\date{}
\date{This work was supported in part by a grant from the University
Research Council of DePaul University.}
\date{}
\maketitle

\begin{abstract}
We describe the prime and primitive spectra for quantized enveloping
algebras at roots of 1\ in characteristic zero in terms of the prime
spectrum of the underlying enveloping algebra. \ Our methods come from the
theory of Hopf algebra crossed products. For primitive ideals we obtain an
analogue of Duflo's Theorem, which says that every primitive ideal is the
annihilator of a simple highest weight module. \ This depends on an
extension of Lusztig's tensor product theorem.
\end{abstract}

\section{Introduction}

We describe the prime and primitive spectra of the restricted specialization
of the quantized enveloping algebra $U_{\zeta }$ where $\zeta $ is a root of
unity of odd order $\ell $ (prime to $3$ in type $G_{2})$ associated to a
semisimple complex Lie algebra $\mathfrak{g}$ of rank $n$ as in e.g. [Lu1,
Lu2, Lu4, Li, An], \ working over a field $\mathrm{k}$ of characteristic
zero \ We refer to $U_{\zeta }$ as the hyperalgebra because it is defined
using a divided power form in a way which mimics the usual hyperalgebra
construction for an affine algebraic group.

There is a Hopf algebra map $\limfunc{Fr}:U_{\zeta}\rightarrow U$ where $U$
is the ordinary enveloping algebra, known as the quantum Frobenius map. \ It
is known [Li, An] that the Hopf kernel of $\limfunc{Fr}$ is the quantized
restricted enveloping algebra $u_{\zeta}$, which is a normal sub- Hopf
algebra of $U_{\zeta}$. Furthermore $U_{\zeta}$ is isomorphic as an algebra
to the crossed product $u_{\zeta}\#_{t}U$ as $\mathrm{k}$-algebras for some
2-cocycle $t:U\otimes U\rightarrow u_{\zeta}$ \ This is explicitly shown in
[An], follows from the general results in [Sch], and was pointed out
informally by Takeuchi earlier. See [Mo, Do] for details about crossed
products. \ The crossed product decomposition is an important tool for our
study of prime and primitive ideals.

The crossed product structure allows to describe the prime spectrum of $%
U_{\zeta }$ in terms of the prime spectra of $u_{\zeta }$ and $U$. \ The
structure of the latter as a topological space is only partially understood
[So], though it is explicitly described in rank one [NG] and in type $A_{2}$
[Ca]. \ We think of the crossed product structure as being a ring-theoretic
explanation of the Lusztig twisted tensor product theorem [Lu1, Li]. Indeed,
we use it to give a different proof of the tensor product theorem which
moreover applies to all simple modules. This is then applied to prove a
version of Duflo's Theorem, which states that every primitive ideal is the
annihilator of a simple highest weight module. \ This result invites the
study of primitive spectra for $U_{\zeta }$\thinspace in connection with
Weyl group orbit geometry, as has been carried out extensively for $U$, see
[Ja] (and [Jo] 8.5.4 for more references). \ We note that the generic
quantized version of Duflo's Theorem was proved by Joseph and Letzter, see
[Jo].

In the second section we discuss $U_{\zeta }$ and its images modulo
annihilators of a simple $u_{\zeta }$-modules. The crossed product in the
image is shown to be a tensor product with trivial action and cocycle, facts
which depend on the Skolem-Noether Theorem and the vanishing of the Sweedler
two-cocycle. We use this fact to show that the prime spectrum of $U_{\zeta }$
is homeomorphic to $\limfunc{Spec}(U)\times \limfunc{Spec}(u_{\zeta .})$ in
section 3. \ We turn to primitive ideals in section 4, first developing the
theory of weights, and simple highest weight modules. Here we allow for
possibly infinite-dimensional simple modules, which are labeled by algebra
maps from the Cartan part $U_{\zeta }^{0}$ to the base field k. The
resulting group of weights has a very different structure from the generic
case. In passing, we give an expression for the primitive elements in $%
U_{\zeta }$, which is not used in what follows, but may be of independent
interest since they are (in a Hopf-algebraic sense) a natural choice as
generators. In section 4.2, assuming k algebraically closed, we prove an
extension of Lusztig's tensor product theorem [Lu1,CP, Li] to all simple
modules for $U_{\zeta }.$ We provide the version of Duflo's theorem in the
last subsection, 4.3.

The authors would like to thank I. M. Musson for his suggestions and Z. Lin
for useful discussions.

\subsection{ Notation and Preliminaries}

We adopt the notation and conventions for the quantized enveloping algebra $%
U_{\zeta }$ of rank $n$ at the primitive root of unity $\zeta $ of order $%
\ell $ as in [Lu1,2,4] (also [An, Li]), working over a field $\mathrm{k}$ of
characteristic zero. $\ell $ is odd and not divisible by $3$ if there is a
component of type $G_{2}.$ The algebra $U_{\zeta }$ has generators $%
E_{i},F_{i},E_{i}^{(\ell )},F_{i}^{(\ell )},K_{i},\binom{K_{i}}{\ell }$, $%
i=1,...n.$ Let $(a_{ij})$ be the symmetrizable Cartan matrix of finite type
with integers $d_{i}\in \{1,2,3\}$ such that $d_{i}a_{ij}=d_{j}a_{ji}$ for
all $i,j.$ Let $\zeta _{i}=\zeta ^{d_{i}}.$ We invoke the relations $%
K_{i}^{\ell }=1$ (as done in e.g. [Li]), and thus only consider type 1
modules. \ We shall use $\Delta $ and $\varepsilon $ for comultiplication
and counit, using subscription to indicate the Hopf algebra in question.

Let $L(m)$ denote the finite-dimensional simple highest weight $U_{\zeta }$%
-module with restricted integral weight $m\in \lbrack 0,\ell -1]^{n}.$ Then
the restriction of this module to $u_{\zeta }$ is simple, and every simple $%
u_{\zeta }$-module is of this form. The proof of these facts can be found in
[Lu1, Lu2] (cf. [CP, 11.2.10]) in the simply-laced case, and it is
well-known that they hold in general.

Working over $\mathbb{Z}[q,q^{-1}],$ let $q_{i}=q^{d_{i}}.$ For $c,m\in 
\mathbb{Z}$ and $r,t\in \mathbb{Z}^{+}$ set 
\begin{eqnarray*}
\QATOPD[ ] {m}{t}_{q_{i}} &=&\prod\limits_{s=1}^{t}\frac{%
q_{i}^{m-s+1}-q_{i}^{-(m-s+1)}}{q_{i}^{s}-q_{i}^{-s}} \\
\QATOPD[ ] {K_{i};c}{t}_{q_{i}} &=&\prod\limits_{s=1}^{t}\frac{%
K_{i}q_{i}^{c-s+1}-K_{i}^{-1}q_{i}^{-(c-s+1)}}{q_{i}^{s}-q_{i}^{-s}}.
\end{eqnarray*}

Specializing $q$ to $\zeta $, we have the corresponding elements $\QATOPD[ ]
{m}{t}_{\zeta _{i}}\in $k and $\QATOPD[ ] {K_{i};c}{t}_{\zeta _{i}}$ $\in $ $%
U_{\zeta }$. Thus $\binom{K_{i}}{\ell }=\QATOPD[ ] {K_{i};0}{\ell }_{\zeta
_{i}}\in U_{\zeta }.$

We shall need several formulas for later use that we summarize here. Fix $%
1\leq i,j\leq n$. Then we have the defining conjugation relations 
\begin{eqnarray}
K_{i}E_{j}K_{i}^{-1} &=&\zeta _{i}^{a_{ij}}E_{j}  \notag \\
K_{i}F_{j}K_{i}^{-1} &=&\zeta _{i}^{-a_{ij}}F_{j}.  \notag
\end{eqnarray}%
These relations imply

\begin{eqnarray}
\binom{K_{i}}{\ell }F_{j}^{(t)} &=&F_{j}^{(t)}\binom{K_{i};-ta_{ij}}{\ell }
\label{3} \\
\binom{K_{i}}{\ell }E_{j}^{(t)} &=&E_{j}^{(t)}\binom{K_{i};ta_{ij}}{\ell }.
\label{4}
\end{eqnarray}%
The next equation is a simple verification:%
\begin{equation}
\QATOPD[ ] {m}{t}_{q_{i}}=(-1)^{t}\QATOPD[ ] {-m+t-1}{t}_{q_{i}}.  \label{5}
\end{equation}

For $m\in \mathbb{Z}$, we write $m=m_{0}+m_{1}\ell \in \mathbb{Z}$ with
unique $m_{0},m_{1}\in \mathbb{Z}$ with $0\leq m_{0}\leq \ell -1$ and refer
to this as the "short $\ell $-adic decomposition of $m$". We mention the
well-known fact that%
\begin{equation}
\QATOPD[ ] {m}{\ell }_{\zeta _{i}}=m_{1}\text{,}  \label{6}
\end{equation}%
see [Lu1], also [CP, 9.2].

Let $m$ be an integer with $0\leq m\leq \ell -1$, and let $c=c_{0}+c_{1}\ell 
$ be the short $\ell $-adic decomposition of $c$. Then we have%
\begin{eqnarray}
\QATOPD[ ] {m-c}{\ell }_{\zeta _{i}} &=&\QATOPD\{ . {-c_{1}\text{ if }m\geq
c_{0}}{-(c_{1}+1)\text{ if }m<c_{0}}  \label{7} \\
\QATOPD[ ] {m+c}{\ell }_{\zeta _{i}} &=&\QATOPD\{ . {c_{1}\text{ if }%
m+c_{0}<\ell }{c_{1}+1\text{ if }m+c_{0}\geq \ell }.  \label{8}
\end{eqnarray}%
The first equation \ref{7} of this last pair of equations can seen by
applying equation \ref{5} to obtain $\QATOPD[ ] {m-c}{\ell }_{\zeta
_{i}}=(-1)^{\ell }\QATOPD[ ] {c_{1}\ell -(m-c_{0})+\ell -1}{\ell }_{\zeta
_{i}}$. If $m-c_{0}\geq 0$, we get $(-1)^{\ell }c_{1}=-c_{1}$ using \ref{6}$.
$ If $m-c_{0}<0$, then $c_{1}\ell -(m-c_{0})+\ell -1=(c_{1}+1)\ell +c_{0}-m-1
$ and we get $-(c_{1}+1)$ for the second case. The second equation \ref{8}
is similar.

\section{Crossed Products}

We begin by briefly review some results concerning cleft extensions and
crossed products. Readers may wish to refer to [Mo, Do] for more detailed
discussions. Fix a subalgebra $A$ of a Hopf algebra $U$ and a right $U$%
-comodule algebra $T$. An extension of algebras $A\subset T$ is said to be 
\textit{cleft}, or more precisely $U$-\textit{cleft,} if there exists a
convolution-invertible $U$-comodule map $\chi \in \mathrm{Hom}(U,T)$ and the
subalgebra of coinvariants $T^{\mathrm{co}U}$ for $U$ equals $A$. Such a map 
$\chi $ is called a \textit{section}. A result from [DT] is that $T$ is
isomorphic to a crossed product $A\#_{t}U$ for some action of $U$ on $A$ and
invertible cocycle $t$ if and only if $T$ is $U$-cleft.

In [Do], it is shown that for a Hopf algebra $U$ and any $U$-cleft extension 
$T$ of the algebra $A$ any two sections $\chi ^{\prime }$, $\chi \in \mathrm{%
Hom}(U,T)$ are related by $\chi ^{\prime }=v\ast $ $\chi $ for some
convolution invertible $v\in \mathrm{Hom}(U,A).$ Recall also that
corresponding to each section $\chi ,$ there is an $\mathrm{k}$-isomorphism 
\begin{equation*}
\tilde{\chi}:A\otimes U\rightarrow T
\end{equation*}%
given by $\tilde{\chi}(a\otimes y)=a\chi (y)$, $a\in A$, $y\in U,$ which is
an algebra isomorphism when the tensor product is given the crossed product
structure. The crossed product is denoted by $A\#_{t}U$ (with action and
twisting $t$ depending on $\chi $) and the elements $a\otimes y$ are written 
$a\#y$ or $a\#_{t}y$. Using the second section, we affix $^{\prime }$ and
obtain an equivalent crossed product and isomorphism $\tilde{\chi}^{\prime
}:A\#_{t^{\prime }}U\rightarrow T$ with new twisting $t^{\prime }$ and 
\begin{equation*}
\tilde{\chi}^{\prime }(a\#_{t^{\prime }}y)=a\chi ^{\prime }(y)=\sum
av(y_{1})\chi (y_{2})
\end{equation*}%
with $a\#_{t^{\prime }}y\in A\#_{t^{\prime }}U.$ The corresponding algebra
isomorphism $A\#_{t^{\prime }}U\rightarrow A\#_{t}U$ sends $a\#_{t^{\prime
}}y$ to $\sum av(y_{1})\#_{t}y_{2}$ [Do, p. 3064].

We now look at $U_{\zeta }\ $and fix notation. The Hopf algebra extension
(see [An]) 
\begin{equation*}
\mathrm{k}\rightarrow u_{\zeta }\rightarrow U_{\zeta }\overset{\limfunc{Fr}}{%
\rightarrow }U\rightarrow \mathrm{k}
\end{equation*}%
is cleft. Here $U=U(\mathfrak{g}),$ $U_{\zeta }$ is a right $U$-comodule
using (id$\otimes \limfunc{Fr})\Delta $ and $u_{\zeta }$ is the coinvariant
subalgebra. The fact that $u_{\zeta }$ is the set of right coinvariants
follows directlty from the PBW theorem. We note for later use that $U_{\zeta
}$ is also a left $U$-comodule using ($\limfunc{Fr}\otimes id)\Delta $ and $%
u_{\zeta }$ is the left coinvariant subalgebra. An explicit section $\gamma
:U\rightarrow U_{\zeta }$ is given by first specifying 
\begin{equation*}
\gamma (e_{i})=E_{i}^{(\ell )},\gamma (f_{i})=F_{i}^{(\ell )},\gamma (h_{i})=%
\binom{K_{i}}{\ell }
\end{equation*}%
where the $e$'s, $f$'s and $h$'s are the usual generators for the enveloping
algebra [An]; we then extend to obtain three algebra maps $U^{+}\rightarrow
U_{\zeta }^{+}$, $U^{-}\rightarrow U_{\zeta }^{-}$, and $U^{0}\rightarrow
U_{\zeta }^{0}$. The definition of $\gamma $ is completed by tensoring,
using the triangular decompositions 
\begin{eqnarray*}
U &=&U^{-}\otimes U^{0}\otimes U^{+} \\
U_{\zeta } &=&U_{\zeta }^{-}\otimes U_{\zeta }^{0}\otimes U_{\zeta }^{+},
\end{eqnarray*}%
i.e. $\gamma (x^{-}x^{0}x^{+})=$ $\gamma (x^{-})\gamma (x^{0})\gamma (x^{+})$
for $x^{+}\in U^{+}$, $x^{0}\in U^{0}$, $x^{-}\in U^{-}.$ It is very
straightforward check that $\gamma $ is both a right and left $U$-comodule
map. The section $\gamma $ corresponds to the crossed product $u_{\zeta
}\#_{t}U\cong U_{\zeta }$ for some action of $U$ on $u_{\zeta }$ and some
2-cocycle $t:U\otimes U\rightarrow u_{\zeta }.$

We remark that it follows from results in [Sch] that a section $U\rightarrow
U_{\zeta }$ exists. Also, our section $\gamma $ differs slightly form the
one in [An], since there the relations $K_{i}^{\ell }=1$ are not assumed.

Let $A$ be an algebra with crossed product action by a Hopf algebra $U.$
Then $T=A\#_{t}U$ is an associative algebra. We briefly review some
ring-theoretic facts and definitions. If $I$ is an $U$-invariant ideal of $A$%
, then is is standard and easy to check that $TJ\subset JT=J\#T$ is an ideal
of $T$. For an ideal $I$ of $A$, as usual (e.g. [Ch]), we put $(I:U)=\{a\in
A|x.a\in I,$ all $x\in U\},$ the largest $U$-invariant ideal of $A$
contained in $I.$ A $U$-invariant ideal $J$ of $A$ is said to be $U$- prime
if $KL\subset J$ implies $K\subset J$ or $L\subset J$, for all $U$-invariant
ideals $L,K\subset A.$

The action of $U$ on $A$ is said to be \textit{inner, }see [Mo or Do] if
there exists a convolution-invertible $u\in $Hom$_{\text{k}}(U,A)$ such that 
$y.a=\sum u(y_{1})au^{-1}(y)$ for all $a\in A$ and $y\in U.$ Here the action
of $U$ is said to be \textit{implemented} by $u$. Note that $u$ is not, in
general, an algebra map.

\begin{lemma}
Consider the crossed product $U_{\zeta }=u_{\zeta }\#_{t}U$ as above. \ Then
the $U$-prime ideals and the prime ideals of $u_{\zeta }$ coincide.
\end{lemma}

\begin{proof}
Let $\mathfrak{m}$ be a $U-$prime ideal of $u_{\zeta }$. \ Then by [Ch], $%
\mathfrak{m}$ is a prime ideal. Conversely let $\mathfrak{m}$ be a prime
ideal of $u_{\zeta }$. \ Then $(\mathfrak{m}:U)$ is a $U$-prime ideal. By
[Ch] again, $(\mathfrak{m}:U)$ is actually a prime and hence maximal ideal.
So $\mathfrak{m}=(\mathfrak{m}:U)$. This shows that every prime ideal of $%
u_{\zeta }$ is $U-$prime.
\end{proof}

\begin{proposition}
\label{crossed product}Consider the crossed product $U_{\zeta }=u_{\zeta
}\#_{t}U$ arising from the section $\gamma $ as above. \ Let \emph{m} be a
prime ideal of $u_{\zeta }$ and set $R$=$u_{\zeta }/\mathfrak{m}$. \ Let $T=$
$U_{\zeta }/\mathfrak{m}U_{\zeta }$ and let $\pi :U_{\zeta }\rightarrow $ $T$
denote the natural map. Then\newline
(a) $\chi =\pi \gamma \in Hom(U,T)$ is a section which yields the crossed
product 
\begin{equation*}
T\cong R\#_{\bar{t}}U
\end{equation*}%
with cocycle $\bar{t}=\pi t:U\otimes U\rightarrow R,$ via an isomorphism $%
\tilde{\chi}:R\#_{t}U\rightarrow T$ arising from the section $\chi $\newline
(b) there is another section $\chi ^{\prime }\in Hom_{k}(U,T)$ giving rise
to the trivial crossed product 
\begin{equation*}
T\cong R\otimes U
\end{equation*}%
via an isomorphism $\tilde{\chi}^{\prime }:R\otimes U\rightarrow T$ of cleft
extensions of $R$.\newline
(c) The corresponding isomorphism $\theta ^{\prime }:R\#_{t}U\tilde{%
\rightarrow}R\otimes U$ is given by 
\begin{equation*}
\theta ^{\prime }(a\#y)=\sum av(y_{(1)})\otimes y_{(2)}
\end{equation*}%
for all $a\in R,$ $y\in U$, for some $v\in Hom_{\mathrm{k}}(U,R)$
\end{proposition}

\begin{proof}
By the Lemma $\mathfrak{m}$ is a $U$-invariant ideal. Therefore, as pointed
out before Lemma 1, $U_{\zeta }\mathfrak{m}\subseteq \mathfrak{m}U_{\zeta }=$
$\mathfrak{m}\#_{t}U,$ so $\mathfrak{m}U_{\zeta }$ is an ideal of $U_{\zeta
}.$ Since $\mathfrak{m}$ is a trivial right $U$ -comodule, $\mathfrak{m}%
U_{\zeta }$ is also a right $U$-subcomodule of $U_{\zeta }.$ Hence $T$
inherits a $U-$comodule algebra structure.

(a) Plainly, $\chi $ is a section of $T$. Thus we have $T\cong R\#_{\bar{t}%
}U $ with cocycle $\bar{t}=\pi t:U\otimes U\rightarrow R.$

(b) and (c): Obviously $R$ is a simple algebra; in fact $\mathfrak{m}$ is
the annihilator of a finite dimensional simple $u_{\zeta }$-module of
dimension, say $n$. \ Since the simple modules are absolutely simple, $R$ is
isomorphic to the algebra of $n\times n$ matrices over $\mathrm{k}$.

The Skolem-Noether Theorem ([Sw], also see [Mo, 6.2]) asserts that the
action of $U$ on $R$ is inner. \ By [Mo, 7.3] we can perform a change of
basis to obtain $R\#_{t}U\cong R_{\tau }[U]$, a twisted product, with
trivial action and a new cocycle $\tau .$ The isomorphism is given
explicitly by a map that is the identity on $R=R\#1$ as follows 
\begin{equation*}
a\#_{t}y\longmapsto \sum au^{-1}(y_{(1)})\#_{\tau }y_{(2)}
\end{equation*}%
for all $y\in U,$ $a\in R,$ where the inner action is implemented by $%
u,u^{-1}\in \limfunc{Hom}(U,R)$ as in [Mo, 7.3].

Furthermore, the action of $U$ is now trivial, so $R$ is a (trivial) $U$%
-module. By [Mo 7.1] or [Do], since $U$ is cocommutative, the cocycle $\tau $
has image in the center of $R$, namely $\mathrm{k}$. The reader may also
check as an exercise that the twisted module condition quickly implies that,
for a trivial action of any Hopf algebra, any cocycle has central values. \
Thus the Sweedler cocycle $\tau \ $corresponds to a Lie cocycle ( [Sw, 4.1
or 4.3]) with coefficients in the trivial module. Therefore $\tau $ is
equivalent to a Sweedler coboundary by Whitehead's second Lemma. \ Here we
have used the equivalence of the Lie (and Hochschild) and Sweedler
cohomologies in degree 2 [Sw, section 4]. We conclude by [Do, 2.4] that \ $%
R_{\tau }[U]$ is isomorphic a smash product with the same trivial action of $%
U$ on $R$; i.e, $R_{\tau }[U]\cong R\otimes U.$ Thus $R\#_{t}U\cong R\otimes
U$ by an isomorphism that is the identity on $R=R\#1=R\otimes 1$.

The isomorphism $R\#_{t}U\cong R\otimes U$ is given by 
\begin{equation*}
a\#y\rightarrow \sum av(y_{(1)})\otimes y_{(2)}
\end{equation*}%
$a\in R,$ $y\in U,$ where $v=u\ast w\in \limfunc{Hom}(U,R)$ and $w$ is the
map whose image is the coboundary $\tau $, by [Do, 2.4]. The corresponding
section $\chi ^{\prime }$ is given by $v\ast \chi $.
\end{proof}

\section{Prime Spectrum}

\begin{theorem}
$\limfunc{Spec}(U_{\zeta})$ is homeomorphic to $\limfunc{Spec}(U)\times%
\limfunc{Spec}(u_{\zeta}).$
\end{theorem}

Proof \ Fix a prime ideal $P\vartriangleleft U_{\zeta }$ and let $\mathfrak{m%
}=P\cap u_{\zeta },$ and put $R=u_{\zeta }/\mathfrak{m}$. \ Then $\mathfrak{m%
}$ is a $U$-prime, and hence maximal, ideal of $R$ by the Lemma. \ As in
[Ch], we define 
\begin{equation*}
\limfunc{Spec}\nolimits_{\mathfrak{m}}(U_{\zeta })=\{Q\in \limfunc{Spec}%
(U_{\zeta }):Q\cap u_{\zeta }=\mathfrak{m}\}.
\end{equation*}%
Notice that the sets $\limfunc{Spec}\nolimits_{\mathfrak{m}}(U_{\zeta })$
are the connected components of $\limfunc{Spec}(U_{\zeta }).$ \ By e.g. [Ch]
we have $\limfunc{Spec}\nolimits_{\mathfrak{m}}(U_{\zeta })\approx \limfunc{%
Spec}(\frac{U_{\zeta }}{\mathfrak{m}U_{\zeta }})\approx \limfunc{Spec}%
(R\#_{t}U)$ \ By the Proposition above we see that $\limfunc{Spec}\nolimits_{%
\mathfrak{m}}(U_{\zeta })\approx \limfunc{Spec}(R\#U)\approx \limfunc{Spec}%
(U $). \ Thus we obtain a mapping $\pi :\limfunc{Spec}(U_{\zeta
})\rightarrow \limfunc{Spec}(U)$ that restricts to a homeomorphism on each $%
\limfunc{Spec}\nolimits_{\mathfrak{m}}(U_{\zeta }).$ \ 

To complete the proof, define the mapping $h:\limfunc{Spec}(U_{\zeta
})\rightarrow\limfunc{Spec}(u_{\zeta})\times\limfunc{Spec}(U)$ by $%
h(P)=(P\cap u_{\zeta})\times\pi(P).$ \ \ Clearly, $h$ maps $\limfunc{Spec}%
\nolimits_{\mathfrak{m}}(U_{\zeta})$ homeomorphically onto the connected
component $\mathfrak{m}\times\limfunc{Spec}(U)$.

\section{Primitive Spectrum}

\subsection{Simple Highest Weight Modules}

\subsubsection{Weights}

Let $U_{\zeta }^{0}$ denote the Hopf subalgebra of $U_{\zeta }$ generated by
the $K_{i}$ and the $\binom{K_{i}}{\ell },$ $i=1,...,n$, which may also be
constructed by specializing the corresponding subalgebra of the Lusztig form
to $\zeta $. \ Also let $G$ denote the multiplicative group of group-likes
generated by the $K_{i}$'s. Since $U_{\zeta }^{0}$ is a pointed commutative
and cocommutative Hopf algebra, by the Kostant structure theorem (cf. [Mo]),
it is the tensor product of the (Hopf) algebras $\mathrm{k}G$ and the
commutative enveloping algebra $U(P(U_{\zeta }^{0}))$, with primitive
generators $P(U_{\zeta }^{0})$). \ On the other hand, restricting the
mapping Fr to $U_{\zeta }^{0}$ and the section $\gamma $ to the subalgebra $%
U^{0}=$ k$[h_{1},...,h_{n}]$, we see that $U_{\zeta }^{0}$ is a crossed
product $\mathrm{k}G\#\mathrm{k}[\binom{K_{1}}{\ell },...,\binom{K_{n}}{\ell 
}]$ with trivial action and trivial cocycle. Thus

\begin{proposition}
$U_{\zeta }^{0}$ is the polynomial algebra generated over $\mathrm{k}G$ by
the algebraically independent elements $\binom{K_{i}}{\ell },$ $i=1,$..$.,n$.
\end{proposition}

The $\binom{K_{i}}{\ell }$ are not the primitive generators of $U_{\zeta
}^{0}$, cf. Remark 1 following Proposition 4 below.

We consider \textit{weights} to be algebra homomorphisms $\lambda :U_{\zeta
}^{0}\rightarrow \mathrm{k}.$ \ We let 
\begin{equation*}
X=Alg_{\mathrm{k}}(U_{\zeta }^{0},\mathrm{k})
\end{equation*}%
denote the set of weights. Weights $\lambda $ are given by a pair of algebra
maps $\lambda =(\lambda _{0},\lambda _{1})$ with $\lambda _{0}:\mathrm{k}%
G\rightarrow \mathrm{k}$ and 
\begin{equation*}
\lambda _{1}:\mathrm{k}[\binom{K_{1}}{\ell },...,\binom{K_{n}}{\ell }%
]\rightarrow \mathrm{k}
\end{equation*}%
and thus are given by a pair of $n$-tuples 
\begin{align*}
(\lambda _{0,1}^{{}},...,\lambda _{0,n}^{{}})& \in \{0,1,...,\ell -1\}^{n} \\
(\lambda _{1,1}^{{}},...,\lambda _{1,n}^{{}})& \in \mathrm{k}^{n}.
\end{align*}%
We say that $\lambda $ is \textit{restricted} weight if $\lambda _{1}=0.$ \
By slight notational abuse we write $\lambda =\lambda _{0}$ in this case,
and similarly $\lambda =\lambda _{1}$ if $\lambda _{0}=\varepsilon _{kG}.$

\bigskip

\textbf{Integral weights}: Let $m$ denote an $n$-tuple of integers $%
(m_{1},...,m_{n})\in \mathbb{Z}^{n}$ We define an embedding $\symbol{126}:%
\mathbb{Z}^{n}\rightarrow X$ of sets by setting 
\begin{eqnarray*}
\tilde{m}(K_{j}) &=&\zeta _{j}^{m_{0j}} \\
\tilde{m}(\binom{K_{j}}{\ell }) &=&m_{1j}^{{}},
\end{eqnarray*}%
where $m_{0j}\ $and $m_{1j}$ are the unique integers in the short $\ell $%
-adic decompositions $m_{0j}+m_{1j}\ell =m_{j}$, $j=1,$..$.,n.$ In our
labeling above of weights as pairs of $n$-tuples, this is equivalent to 
\begin{eqnarray*}
\tilde{m}_{0j} &=&m_{0j} \\
\tilde{m}_{1j} &=&m_{1j}
\end{eqnarray*}%
$j=1,$..$.,n.$

We need to specify the convolution group structure on the set of weights
which arises from the cocommutative Hopf algebra structure on $U_{\zeta
}^{0} $. It is given by \textquotedblleft carrying\textquotedblright\ from
the finite abelian factor to the additive factor as we next state. We use $%
\ast $ for the group operation on $X.$

\begin{proposition}
The group structure on $X=$Alg($U_{\zeta }^{0},\mathrm{k}$) is given by 
\begin{equation*}
((\lambda _{0},\lambda _{1})\ast (\mu _{0},\mu _{1}))_{j}==\QATOPD\{ .
{(\lambda _{0j}+\mu _{0j},\lambda _{1j}+\mu _{1j})\text{ if }\lambda
_{0j}+\mu _{0j}<\ell }{(\lambda _{0j}+\mu _{0j}-\ell ,\lambda _{1j}+\mu
_{1j}+1)\text{ if }\lambda _{0j}+\mu _{0j}\geqslant \ell }
\end{equation*}%
for $\lambda ,$ $\mu \in X;$ $j=1,...,n.$
\end{proposition}

\begin{proof}
Let $\lambda ,$ $\mu \in X$. $\func{Si}$nce each $K_{i}$ is group-like, 
\begin{align*}
\lambda \ast \mu (K_{i})& =\lambda (K_{i})\mu (K_{i}) \\
& =\zeta ^{\lambda _{0i}+\mu _{0i}}=\QATOPD\{ . {\zeta ^{\lambda _{0i}+\mu
_{0i}}\text{ if }\lambda _{0i}+\mu _{0i}<\ell }{\zeta ^{\lambda _{0i}+\mu
_{0i}-\ell }\text{ otherwise}}
\end{align*}

Further, recall (e.g. [CP]) that $\Delta (\binom{K_{i}}{\ell }%
)=\sum\limits_{j=0}^{\ell }\binom{K_{i}}{\ell -j}K^{-j}\otimes \binom{K_{i}}{%
j}K^{\ell -j}$, and observe that $\lambda (\binom{K_{i}}{j})=\binom{\lambda
_{0i}}{j}_{\zeta _{i}}$ for all $j<\ell $. We infer that 
\begin{equation*}
\lambda \ast \mu (\binom{K_{i}}{\ell })=\lambda _{1i}+w+\mu _{1i}
\end{equation*}%
where the cross-term is $w=\sum\limits_{j=1}^{\ell -1}\zeta ^{-(\lambda
_{0i}+\mu _{0i})j}\binom{\lambda _{0i}}{\ell -j}_{\zeta _{i}}\binom{\mu _{0i}%
}{j}_{\zeta _{i}}.$ The proof is completed by applying equation \ref{****}
below to deduce that $w=\binom{\lambda _{0i}+\mu _{0i}}{\ell }_{\zeta _{i}}$%
. By \ref{6} we immediately obtain the desired result that $w=0$ or $w=1$ in
cases as in the statement.
\end{proof}

The following consequence is now straightforward. It explains the
relationship between our weights and the integral weights as considered in
[L1, L2], which were given as $n$-tuples of integers.

\begin{corollary}
The map $\symbol{126}:\mathbb{Z}^{n}\rightarrow X$ defined above is an
embedding of abelian groups.
\end{corollary}

From now on we will use additive notation for the group operation on $X$.
For future use we record the following easy facts.

\begin{proposition}
Let $\lambda =(\lambda _{0},\lambda _{1})$ and $\mu =(\mu _{0},\mu _{1})$ be
weights and write $\lambda _{j}=(\lambda _{0j},\lambda _{1j})$, $j=1,...,n.$
The following equations hold 
\begin{eqnarray*}
(-\lambda )_{j} &=&\QATOPD\{ . {(\ell -\lambda _{0j},-\lambda _{1j}-1)\text{
if }\lambda _{0j}>0}{(0,-\lambda _{1j})\text{ if }\lambda _{0j}=0} \\
(\lambda -\mu )_{j} &=&\QATOPD\{ . {(\ell +\lambda _{0j}-\mu _{0j},\lambda
_{1j}-\mu _{1j}-1)\text{ if }\lambda _{0j}<\mu _{0j}}{(\lambda _{0j}-\mu
_{0j},\lambda _{1j}-\mu _{1j})\text{ if }\lambda _{0j}\geq \mu _{0j}}.
\end{eqnarray*}
\end{proposition}

\begin{remark}
We describe the primitive elements of $U_{\zeta }^{0}.$ It suffices to work
in rank 1, so let $U=U_{\zeta }^{0}(sl(2))$ Then up to a scalar multiple, $U$
has a unique primitive element 
\begin{equation*}
\binom{K}{\ell }+\sum_{i=0}^{\ell -1}a_{i}K^{i}
\end{equation*}%
where 
\begin{equation*}
a_{i}=\frac{1}{\ell ^{2}}\sum_{j=0}^{\ell -1}j\zeta ^{-ij},i=0,1,...,\ell -1.
\end{equation*}
\end{remark}

\begin{proof}
The weight group $X=$Alg($U$, \emph{k}) can be written as the set of pairs $%
(\zeta ^{i},a)$, $a\in $\textrm{k}$,$ $i=0,1,...,\ell -1.$ Let $X^{\prime }$
be the group with the same underlying set but with group structure $%
\left\langle \zeta \right\rangle \times \alpha _{\mathrm{k}}$ where $\alpha
_{\mathrm{k}}$ denotes the additive group of \textrm{k}. Notice that the
coordinate Hopf algebra of $X^{\prime }$ is $U^{\prime }=$k$[K,d]$ where $d$
is primitive and $K$ is a grouplike of order $\ell $. Thus Alg($U^{\prime },%
\mathrm{k})$ is isomorphic to $X^{\prime }$ where the entries in the pairs
correspond to evaluations at $K$ and $d$, respectively. It is easy to see
that $X$ and $X^{\prime }$ are isomorphic via $\eta :X\rightarrow
X^{^{\prime }}$ where 
\begin{equation*}
\eta ((\zeta ^{i},a))=(\zeta ^{i},a+\frac{i}{\ell }).
\end{equation*}%
Therefore we have an isomorphism of Hopf algebras $\eta ^{\ast }:U^{\prime
}\rightarrow U.$

We proceed to calculate the primitive element $\eta ^{\ast }(d)$. Let $%
(\zeta ^{i},a)\in X$. Evaluating, we have 
\begin{align*}
& <(\zeta ^{i},a),\eta ^{\ast }(d)> \\
& =<\eta (\zeta ^{i},a),d> \\
& =<(\zeta ^{i},a+\frac{i}{\ell }),d> \\
& =a+\frac{i}{\ell }.
\end{align*}%
On the other hand, since $\binom{K}{\ell }$ is primitive modulo the
coradical k$[K],$ we may write 
\begin{equation*}
\eta ^{\ast }(d)=b\binom{K}{\ell }+\sum_{i=0}^{\ell -1}a_{j}K^{j}
\end{equation*}%
for some $b,a_{j}\in \mathrm{k}.$ Evaluating at $(\zeta ^{i},a)$ yields the
system of $\ell $ equations 
\begin{equation*}
ba+\sum_{j=0}^{\ell -1}a_{i}\zeta ^{ij}=a+\frac{i}{\ell }
\end{equation*}%
$i=0,...,\ell -1$. Setting $a=0,1$ immediately yields $b=1$. So we have the
system of equations 
\begin{equation*}
\sum_{j=0}^{\ell -1}a_{i}\zeta ^{ij}=\frac{i}{\ell },
\end{equation*}%
which can be solved for the $a_{i}$ by inverting the Vandermonde matrix $%
(\zeta ^{ij}).$ Its inverse is (fortunately, as pointed out to us by J.
Angelos) $\frac{1}{\ell }(\zeta ^{-ij})$, as can easily be checked. This
yields the asserted expressions for the coefficients $a_{i}$.

Since $U_{\zeta }^{0}$ has Krull dimension one and the coradical is finite,
the k-space of primitives of $U_{\zeta }^{0}$ must be of dimension one. The
uniqueness follows from the Kostant structure theorem, see [Mo].
\end{proof}

\subsubsection{Highest Weight Modules}

We proceed to develop the theory of highest weight modules, following the
standard strategy. The main difference is that we allow nonintegral weights,
and the weights are an extension of the additive group of $\mathrm{k}$ by a
finite abelian group $G=<K_{1},...,K_{n}>$.

\begin{definition}
Let $V$ denote a $U_{\zeta }$ -module. We say that $v\in V$ is a \textit{%
primitive vector}\ of weight $\lambda \in X$ if 
\begin{align*}
K_{i}.v& =\lambda (K_{i})v \\
\binom{K_{i}}{\ell }v& =\lambda (\binom{K_{i}}{\ell })v
\end{align*}%
and $E_{i}v=E_{i}^{(\ell )}v=0$ for all $j$. We say that $V$ is a \textit{%
highest weight module} if $V$ it is generated by a primitive vector $v.$
\end{definition}

In this case, $V$ is the sum of its weight spaces 
\begin{equation*}
V_{\mu }=\{v\in V|x.v=\mu (x)v,\text{ all }x\in U_{\zeta }^{0}\}
\end{equation*}
for weights\ $\mu .$

We define a partial ordering on $X$ next. Let $(a_{ij})$ be the
symmetrizable Cartan matrix with entries $d_{i}a_{ij}=d_{j}a_{ji}$ in the
standard notation. \textit{Simple Roots} are defined to be the restricted
weights $\rho _{i}=(a_{1i,}a_{2i},...,a_{ni})$, the $i^{th}$ column of $%
(a_{ij}),$ $i=1,..,n$. We say that $\lambda \preceq \mu $ if $\mu -\lambda $
is a nonegative integral linear combination of simple roots. This ordering
is sensible because of the following proposition. Recall that we are using
additive notation for the convolution group operation on $X.$

\begin{proposition}
\label{highest weights}Let $V$ be a $U_{\zeta }$-module and let $v\in V$ be
a vector of weight $\lambda \in X$. Then $F_{i}^{(t)}v$ has weight $\lambda
-t\rho _{i}$ and $E_{i}^{(t)}v$ has weight $\lambda +t\rho _{i}$ for all $%
t\in \mathbb{N}$.
\end{proposition}

\begin{proof}
Observe that 
\begin{align*}
K_{j}F_{i}^{(t)}v& =(K_{j}F_{i}^{(t)}K_{j}^{-1})(K_{j}v) \\
& =(\zeta _{j}^{-ta_{ji}}F_{i}^{(t)})(\zeta _{j}^{\lambda _{0j}}v) \\
& =\zeta _{j}^{\lambda _{0j}-ta_{ji}}F_{i}^{(t)}v
\end{align*}

This establishes the first equation to be checked. We leave the analogous
formula for $K_{j}E_{i}^{(t)}v$ to the reader.

The second pair of equations to be checked are 
\begin{eqnarray*}
\binom{K_{j}}{\ell }F_{i}^{(t)}v &=&(\lambda -t\rho _{i})F_{i}^{(t)}v \\
\binom{K_{j}}{\ell }E_{i}^{(t)}v &=&(\lambda +t\rho _{i})E_{i}^{(t)}v\text{.}
\end{eqnarray*}%
We begin by recalling a pair of identities (e.g. [Lu2, g9-g10], p. 270). Let 
$m,c\in \mathbb{N}.$ Then%
\begin{align}
\QATOPD[ ] {K_{j};-c}{\ell }_{q_{j}}& =\sum\limits_{s=0}^{\ell
}(-1)^{s}q_{j}^{c(\ell -s)}\QATOPD[ ] {c+s-1}{s}_{q_{j}}K_{j}^{s}\QATOPD[ ] {%
K_{j}}{\ell -s}_{q_{j}}  \label{*} \\
\QATOPD[ ] {K_{j};c}{\ell }_{q_{j}}& =\sum\limits_{s=0}^{\ell }q_{j}^{c(\ell
-s)}\QATOPD[ ] {c}{s}_{q_{j}}K_{j}^{-s}\QATOPD[ ] {K_{j}}{\ell -s}_{q_{j}}.
\label{**}
\end{align}%
The next pair of formulas are obtained from these identities by specializing 
$K_{j}$ to $q_{j}^{m}:$ 
\begin{align}
\QATOPD[ ] {m-c}{\ell }_{q_{j}}& =\sum\limits_{s=0}^{\ell
}(-1)^{s}q_{j}^{\ell c+s(m-c)}\QATOPD[ ] {c+s-1}{s}_{q_{j}}\QATOPD[ ] {m}{%
\ell -s}_{q_{j}}  \label{***} \\
\QATOPD[ ] {m+c}{\ell }_{q_{j}}& =\sum\limits_{s=0}^{\ell }q_{j}^{\ell
c-s(m-c)}\QATOPD[ ] {c}{s}_{q_{j}}\QATOPD[ ] {m}{\ell -s}_{q_{j}}.
\label{****}
\end{align}

If we put $m=\lambda _{0j}$, the $s=0$ term in equations \ref{***} and \ref%
{****} is zero. Let $S$ denote the sum of the terms in equations \ref{*} or %
\ref{**} with $s>0.$ Now using the obvious fact that $\binom{K_{j}}{\ell -s}%
v=\QATOPD[ ] {\lambda _{0j}}{\ell -s}_{\zeta _{j}}v$ for all $s>0$, we see
that 
\begin{equation*}
S.v=\QATOPD[ ] {\lambda _{0j}-c}{\ell }_{\zeta _{j}}v.
\end{equation*}%
Therefore 
\begin{align*}
\QATOPD[ ] {K_{j};-c}{\ell }_{\zeta _{j}}v& =\binom{K_{j}}{\ell }v+\QATOPD[ ]
{\lambda _{0j}-c}{\ell }_{\zeta _{j}}v \\
& =(\lambda _{1j}+\QATOPD[ ] {\lambda _{0j}-c}{\ell }_{\zeta _{j}})v.
\end{align*}%
Set $c=ta_{ji}$. Now by equation \ref{3},%
\begin{equation*}
\binom{K_{j}}{\ell }F^{(t)}.v=(\lambda _{1j}+\QATOPD[ ] {\lambda _{0j}-c}{%
\ell }_{\zeta _{j}})F^{(t)}.v.
\end{equation*}%
Writing the short $\ell $-adic decomposition $c=$ $c_{0}+c_{1}\ell $, one
has by \ref{7} 
\begin{equation*}
\QATOPD[ ] {\lambda _{0j}-c}{\ell }_{\zeta _{j}}=\QATOPD\{ . {-c_{1}\text{
if }\lambda _{0j}\geq c_{0}}{-(c_{1}+1)\text{ otherwise}}.
\end{equation*}%
This is the desired result for the action of the $F_{i}^{(t)}$. The
analogous result for the action of $E_{i}^{(t)}$ is obtained similarly using
equations \ref{4}, \ref{8} and \ref{****} in place of \ref{3}, \ref{7} and %
\ref{***}.
\end{proof}

\begin{corollary}
Let $v$ be a primitive vector of weight $\lambda \in X$. Then $\lambda $ is
the highest weight in $U_{\zeta }v$ in the ordering $\preceq $ defined above.
\end{corollary}

\subsubsection{Simple Highest Weight Modules}

Given a weight $\lambda ,$ we can construct a highest weight module for $%
U_{\zeta }$ by the usual Verma module construction. \ Using the standard
argument, we find that there is a unique simple factor module which we
denote by $L(\lambda )$ with primitive vector of weight $\lambda ,$ and
every simple highest weight module is obtained in this way. Let 
\begin{equation*}
T_{\lambda }=U_{\zeta }/\mathrm{ann}_{u_{\zeta }}(L(\lambda ))U_{\zeta }
\end{equation*}%
Clearly the representation $U_{\zeta }\rightarrow \mathrm{End}_{\mathrm{k}%
}(L(\lambda ))$ factors through $T_{\lambda }.$ We write 
\begin{equation*}
\phi _{\lambda }:T_{\lambda }\rightarrow \mathrm{End}_{\mathrm{k}}(L(\lambda
))
\end{equation*}%
for the representation of $T_{\lambda }.$

Consider a restricted weight $\lambda $; then $L(\lambda )$ is simple as a $%
u_{\zeta }$-module. Let 
\begin{equation*}
R_{\lambda }=u_{\zeta }/\mathrm{ann}_{u_{\zeta }}(L(\lambda ))\cong \mathrm{%
End}_{\mathrm{k}}(L(\lambda ))
\end{equation*}

On the other hand, when $\lambda=\lambda_{1}$, then $\lambda$ factors
through $U^{0}=\mathrm{k}[h_{1},...,h_{n}]$ and thus $\lambda$ corresponds
in the obvious way to a weight $\bar{\lambda}$ for $U$. (This makes sense
since $\lambda$ vanishes on $\ker\limfunc{Fr}\cap U_{\zeta}^{0}=\tsum
(K_{i}-1)U_{\zeta}^{0}.)$ The corresponding simple $U$-module with highest
weight $\bar{\lambda}=\bar{\lambda}_{1}\ $shall be denoted by $L_{U}(\bar{%
\lambda})$, so that $L_{U}(\bar{\lambda})^{\limfunc{Fr}}\cong L(\lambda).$

We will often drop the subscript $\lambda$.

\subsection{Tensor Product Theorem for Simple Modules}

We recall the algebra map $\chi ^{\prime }:U\rightarrow T$ from Proposition
1.

\begin{lemma}
Let $\lambda$ be a restricted weight. \ Then $U$ acts trivially on $%
L(\lambda)$ via $\phi_{\lambda}\chi^{\prime}:U\rightarrow\mathrm{End}_{%
\mathrm{k}}(L(\lambda)),$ i.e., $\phi_{\lambda}\chi^{\prime}=\varepsilon
_{U}.$
\end{lemma}

\begin{proof}
By Proposition 1, the elements of $\phi_{\lambda}\chi^{\prime}(U)$ commute
with $R$ (as endomorphisms of $L(\lambda)$). Therefore $\phi_{\lambda}\chi^{%
\prime}$ sends $U$ to $\mathrm{k=End}_{R}(L(\lambda))\subset $\textrm{End}$_{%
\mathrm{k}}$($L(\lambda)$).
\end{proof}

Let $R$ and $U$ be rings, and let $L$ and $V$ be $R$ and $U$ -modules,
respectively. We give $L\otimes V$ (resp. $V\otimes L$) the structure of an $%
R\otimes U$ -module (resp. $U\otimes R$) where the action is along tensor
factors. If $L$ is a $U_{\zeta }-$module and $V$ is a $U$-module, we give $%
L\otimes V^{\limfunc{Fr}}$ and $V^{\limfunc{Fr}}\otimes L$ the usual module
structure via the comultiplication of $U_{\zeta }$, with $V^{\limfunc{Fr}}$
considered as \ a $U_{\zeta }$-module along $\limfunc{Fr}$. \ When $%
R=u_{\zeta }/\mathrm{ann}_{u_{\zeta }}(L)$ we have both module structures on
the vector space $L\otimes V$, because $L$ can be viewed as a $U_{\zeta }$%
-module.

\bigskip

We shall assume k is algebraically closed for the remainder of the paper

\begin{theorem}
(a) Every simple $U_{\zeta }$- module is isomorphic to $L\otimes V^{\limfunc{%
Fr}}$ for some simple highest weight module $L\cong L(\lambda _{0})$ with
restricted weight $\lambda _{0}$ and a simple $U$-module $V$ whose $U_{\zeta
}$-module structure is via $\limfunc{Fr}.$ \newline
(b) Conversely, every module of this form, with $L$ and $V$ simple, is a
simple $U_{\zeta }$- module.\newline
(c) In addition, $V^{\limfunc{Fr}}\otimes L\cong L\otimes V^{\limfunc{Fr}%
}\cong L\otimes V\cong V\otimes L$.
\end{theorem}

\begin{proof}
Let $M$ be a simple $U_{\zeta}$- module. \ Then, as the annihilator $%
\mathfrak{m}=\mathrm{ann}_{u_{\zeta}}(M)$ is a prime ideal, $M$ is a
semisimple $u_{\zeta}$ -module with unique isotypic component, say $%
L=L(\lambda_{0})$.

Let $\mathrm{e}$ be a primitive idempotent in $R=u_{\zeta }/\mathfrak{m}$.
Then $L\cong R\mathrm{e}$ and $\mathrm{e}R\mathrm{e}=\mathrm{k}.$ \ As in
Proposition 1 we have $T=U_{\zeta }/\mathfrak{m}U_{\zeta }\cong R\otimes
U\cong M_{n}(U).$ \ Accordingly, $\mathrm{e}T\mathrm{e}\cong U$ and $T%
\mathrm{e}$ is a projective indecomposable module for $T$. Clearly $L$ is a
generator for mod($R$), so $T\mathrm{e}$ is a generator for module category
Mod($T$).

By Morita equivalence, as left $T$ -modules, 
\begin{align*}
M & \cong T\mathrm{e}\otimes_{U}\mathrm{e}M \\
& \cong(L\otimes U)\otimes_{U}\mathrm{e}M \\
& \cong L\otimes\mathrm{e}M
\end{align*}
where, in the last expression, $T=R\otimes U$ acts along tensor factors, and 
$\mathrm{e}M=V$ is the corresponding simple $\mathrm{e}T\mathrm{e}\cong U$%
-module.

Consider $L\otimes V$ as a $U_{\zeta }$-module via $\pi :U_{\zeta
}\rightarrow T=R\otimes U$, and $L\otimes V^{\limfunc{Fr}}$ as a $U_{\zeta }$%
-module by comultiplying. Note that $L\otimes V^{\limfunc{Fr}}$ inherits a $%
T $-module structure since $\mathfrak{m}U_{\zeta }$ is a right $U$%
-subcomodule. We finish the proof of (a) by showing that\ $L\otimes V\cong
L\otimes V^{\limfunc{Fr}}$. \ To do this it suffices to show that $\chi
^{\prime }(U)$ only acts on the right factor $V^{\limfunc{Fr}}$, i.e. 
\begin{equation*}
\chi ^{\prime }(a)(l\otimes v^{\limfunc{Fr}})=l\otimes av
\end{equation*}
\ for all $a\in U,$ $l\in L$ and $v\in V.$ \ \ To show this, recall that $%
U_{\zeta }$ is right $U$-comodule via $\rho _{\zeta }=(1\otimes \limfunc{Fr}%
)\Delta _{U_{\zeta }}.$ Thus $T$ inherits a right $U$ -comodule with
structure map denoted by $\rho :T\rightarrow T\otimes U$. We summarize these
maps in the following diagram: 
\begin{equation*}
\begin{array}{ccc}
U_{\zeta } & \overset{\Delta }{\longrightarrow } & U_{\zeta }\otimes
U_{\zeta } \\ 
\downarrow _{id} &  & \downarrow _{\pi \otimes \limfunc{Fr}} \\ 
U_{\zeta } & \overset{\rho _{\zeta }}{\longrightarrow } & U_{\zeta }\otimes U
\\ 
\downarrow _{\pi } &  & \downarrow _{\pi \otimes id_{U}} \\ 
T & \overset{\rho }{\rightarrow } & T\otimes U.%
\end{array}%
\end{equation*}%
Similarly $T$ is a left $U$-comodule with structure map $\lambda
:T\rightarrow U\otimes T$.

We now compute $\chi ^{\prime }(a)(l\otimes v^{\limfunc{Fr}})$%
\begin{align*}
& =(\rho \chi ^{\prime }(a))(l\otimes v) \\
& =(\chi ^{\prime }\otimes id_{T})\Delta _{U}(a)(l\otimes v) \\
& =(\varepsilon _{U}\otimes id_{T})\Delta _{U}(a)(l\otimes v) \\
& =(id_{U}\otimes a)(l\otimes v)
\end{align*}%
where the second equality holds because $\chi ^{\prime }$ is a right $U$%
-comodule map, and the third equality holds by the preceding Lemma. This
completes the proof of (a).

To prove (b), let $V$ be a simple $U$-module and $L$ a simple $u_{\zeta }$%
-module. \ By a well-known result [Qu] (since k is algebraically closed), $L$
and $V$ have trivial endomorphism rings $\mathrm{k}$. \ This fact, with a
straightforward density theorem argument, shows that $L\otimes V$ is a
simple $T$-module, with $R$ and $U$ acting along factors. Thus, as a $%
U_{\zeta }$ -module via $\pi :U_{\zeta }\rightarrow T$, we see that $%
L\otimes V$ is a simple module.

To prove (c) we first claim that the right section $\chi ^{\prime }$ is also
a left $U$-comodule map. As in Section 2, we have $U_{\zeta }=u_{\zeta
}\gamma (U)$ and $u_{\zeta }$ is both left and right $U$-coinvariant for the
respective left and right coactions $\lambda _{\zeta }=(\limfunc{Fr}\otimes
id_{U_{\zeta }})\Delta _{U_{\zeta }}$ and $\rho _{\zeta }=(id_{U_{\zeta
}}\otimes \limfunc{Fr})\Delta _{U_{\zeta }}$. Let $a\in u_{\zeta }$ and $%
y\in U$. Then in Sweedler notation, 
\begin{equation*}
\rho _{\zeta }(a\gamma (y))=\sum a\gamma (y_{1})\otimes y_{2}
\end{equation*}%
since $a$ is right coinvariant and $\gamma $ is a right $U$-comodule map; on
the other hand, 
\begin{equation*}
\lambda _{\zeta }(a\gamma (y))=\sum y_{1}\otimes a\gamma (y_{2})
\end{equation*}%
since $a$ is left coinvariant and $\gamma $ is a right $U$-comodule map.
Since $U$ is cocommutative, we see that $\lambda _{\zeta }=\tau \rho _{\zeta
}$, where $\tau $ denotes the twist map. The same fact pushes down to $T$,
i.e., $\lambda =\tau \rho $. Thus for any right comodule map $\varkappa
:U\rightarrow T$, 
\begin{eqnarray*}
\lambda \varkappa  &=&\tau \rho \varkappa  \\
&=&\tau (\varkappa \otimes id_{U})\rho  \\
&=&(id_{U}\otimes \varkappa )\lambda .
\end{eqnarray*}%
This establishes the claim that $\varkappa $ is a left $U$-comodule map

By the left-comodule version of the argument just used in the last paragraph
of the proof of part (a), we see that $V\otimes L\cong V^{\limfunc{Fr}%
}\otimes L$. Obviously, $V\otimes L\cong L\otimes V$ as $U_{\zeta }$%
-modules. Thus $V^{\limfunc{Fr}}\otimes L\cong V\otimes L\cong L\otimes
V\cong L\otimes V^{\limfunc{Fr}}$. This completes the proof of (c) and the
Theorem.
\end{proof}

\begin{remark}
The tensor symmetry in (c) follows from the quasi-tensor structure on $%
U_{\zeta }$ (see [Lu3]) when the modules are finite-dimensional. \ We do not
know this to be the case in general.
\end{remark}

\subsection{An Analog of Duflo's Theorem}

\begin{theorem}
Every primitive ideal of $U_{\zeta }$ is the annihilator of a simple highest
weight module.
\end{theorem}

\begin{proof}
Let $M$ be a simple left $U_{\zeta }-$module \ \ Let $P$ denote the
annihilator in $U_{\zeta }$ and set \emph{m}=$u_{\zeta }\cap P,$ which is a
prime ideal of $u_{\zeta }$ by Lemma 1. \ As a $U_{\zeta }-$module $M$ $%
\cong L\otimes V^{\limfunc{Fr}}$ where $L$ and $V$ are simple $u_{\zeta }$-
and $U$-modules respectively. By Duflo's Theorem [Duf], $\mathrm{ann}_{U}(V)$
is the annihilator of a highest weight module for $U$, say $L_{U}(\bar{%
\lambda}_{1})$. Therefore \textrm{ann}$_{U_{\zeta }}(V^{\limfunc{Fr}})=$ $%
\mathrm{ann}_{U_{\zeta }}(L(\lambda _{1}))$ for some weight $\lambda _{1}$
(notation as in 4.1).

Note also that $L$ is a simple highest weight module for $u_{\zeta }$, say
of highest weight $\lambda _{0}.$ Thus if $x\in L=L(\lambda _{0})$ \ and $%
y\in L_{U}(\bar{\lambda}_{1})$ are primitive vectors for the respective $%
u_{\zeta }$- and $U$-module structures, then adopting the argument in [Lu1], 
\begin{align*}
E_{i}.(x\otimes y)& =(1\otimes E_{i}+E_{i}\otimes K).(x\otimes
y)=E_{i}.x\otimes y=0 \\
E_{i}^{(\ell )}.(x\otimes y)& =(1\otimes E_{i}^{(\ell )}+E_{i}^{(\ell
)}\otimes 1).(x\otimes y)=E_{i}^{(\ell )}.x\otimes y+1+e_{i}.y=0 \\
& \text{and} \\
K_{i}.(x\otimes y)& =\lambda _{0}(K_{i})(x\otimes y)=\lambda _{0,i}(x\otimes
y) \\
\binom{K_{i}}{\ell }.(x\otimes y)& =\lambda _{1}(\binom{K_{i}}{\ell }%
)(x\otimes y)=\lambda _{1,i}(x\otimes y)
\end{align*}%
for all $i$. Hence $x\otimes y\in L\otimes L(\lambda _{1})^{\limfunc{Fr}}$
is a primitive vector, and thus $L\otimes L_{U}(\bar{\lambda}_{1})^{\limfunc{%
Fr}}$ is a simple highest weight module. Finally \textrm{ann}$_{U_{\zeta
}}(M)$ $=$ann$(L(\lambda _{0})\otimes V^{\limfunc{Fr}})$ can be computed
using the Sweedler wedge in the coalgebra $U_{\zeta }$ (see [Mo]) to be 
\begin{align*}
& \text{ann}_{U_{\zeta }}(L(\lambda _{0}))\wedge \text{ann}_{U_{\zeta }}(V^{%
\limfunc{Fr}})) \\
& =\text{ann}_{U_{\zeta }}(L(\lambda _{0}))\wedge \text{ann}_{U_{\zeta
}}(L_{U}(\bar{\lambda}_{1})^{\limfunc{Fr}}) \\
& =\text{ann}_{U_{\zeta }}(L(\lambda _{0})\otimes L_{U}(\bar{\lambda}_{1})^{%
\limfunc{Fr}}))
\end{align*}%
This completes the proof of the Theorem.
\end{proof}

\end{document}